\documentclass[letterpaper]{article}
\usepackage{amssymb}
\usepackage{amsmath}	
\usepackage{mathtools}
\usepackage[utf8]{inputenc}
\usepackage[english]{babel}
\usepackage{amsthm}
\usepackage{tikz,pgf,pgfplots}
\usepackage[utf8]{inputenc}
\usepackage{amsmath}
\usepackage{algorithm}
\usepackage[noend]{algpseudocode}
\usepackage{listings}
\usepackage{courier}

\newtheorem{remark}{Remark}
\newtheorem{definition}{Definition}

\begin{document}

\title{Wirtinger's Calculus for the Load Flow in Power Distribution Grids}
\author{Alejandro Garc\'es \\
ORCiD:0000-0001-6496-0594\\
Universidad Tecnologica de Pereira,\\
post code 660003,\\
alejandro.garces@utp.edu.co
}

\maketitle
\begin{abstract}
This short paper presents a Wirtinger's-Calculus based load-flow methodology for power distribution grids.  This approach allows to obtain an algorithm which works directly on the complex domain maintaining some useful symmetries and a compact representation.  The paper aims to introduce Wirtinger's Calculus as a suitable tool for power systems analysis and simulations; therefore, it is presented as a tutorial, playing especial attention on the implementation in modern scrip languages such as Matlab/Octave, which allow easy representation of complex-arrays with fast calculations.  A Newton's-based method is proposed in which the Jacobian is replaced by Wirtinger's derivatives obtaining a compact representation.  Simulation results complement the analysis. Despite being a mature theory, Wirtinger's-Calculus has not been applied before in this type of problems.
\end{abstract}

\section{Introduction}

Load-flow calculations are indispensable in power systems operation, planning and control. Early algorithms were based on the Newton's and Gauss methods, easily implementable in digital computers.  These methods were modified in order to obtain fast and efficient algorithms such as the decoupled and fast-decoupled load-flow \cite{elhawary}.  However, these approximations showed to be inadequate in power distribution grids for which new methodologies were developed such as the backward/forward sweep load flow\cite{renato}, perhaps the most popular algorithm among them. 

Despite being a classical research area, there is a renewed interest amongst the scientific community, in studding new aspects of the load-flow problem in power distribution grids. Recent investigations have demonstrated the possibility to obtain linear and convex approximations to the power-flow in distribution grids with low demand.  Indeed, at least three linear approximations in power distribution grids have been recently suggested standing out the one proposed by Bolognani \cite{Bolognani} and Marti \cite{marti} as well as the methodology proposed on \cite{yomismo} which is based on a Taylor expansion on the complex plane.  Convex analysis for the load-flow has been also an active research area in which semi-definite and second order cone optimization are the most prominent examples \cite{low}.   Formal analysis about the existence of the solution and convergence of the numerical methods have been also presented; in this aspect, we can highlight the works proposed in \cite{unicidad1}\cite{unicidad2}\cite{unicidad3}.

On the other hand, most of the methodologies in ac power distribution grids are based in real analysis despite being a problem on a complex domain $\mathbb{C}^n$.  In that conventional approach, the power-flow equations are divided into real and imaginary part obtaining a model on a Euclidean domain with a double dimensionality $\mathbb{R}^{2n}$, and losing algebraic properties of the complex domain. Complex analysis would be of course, the most natural approach for power-flow problems; however, the power-flow equations are non-analytic functions and hence we cannot obtain a direct Taylor series expansion on the complex domain. In this context Wirtinger’s calculus could provide an alternative formulation.

A general theory for non-analytical complex functions comes back to the initial works of Poincar\'e further developed by Wilhelm Wirtinger\cite{Wirtinger1927}. This theory deals with non-analytical functions in the complex domain by defining new operators named Wirtinger’s derivatives.  These operators do not fulfill all the properties of conventional complex derivative (for example the Cauchy-Riemann conditions);  This is perhaps the main reason why Wirtinger’s calculus is few known even in the mathematical community.  However, Wirtinger’s calculus could present advantages for power systems analysis at both, high-power and power distribution level.  It could be also important for micro-grids operation. In addition,  Wirtinger’s calculus can simplify the implementation of load-flow algorithms in modern scripting languages such as Matlab/Octave and Python. Recent applications of Wirtinger's calculus outside power systems researches include gradient descent algorithms for machine learning \cite{wit_neural} and signal processing algorithms \cite{signal}.  However, to the best knowledge of the authors, it has not been applied in computer simulations for power systems.

This paper presents a Newton's method for the load-flow in power distribution grids using a Wirtinger's calculus formulation. The proposed formulation is clearly more compact than the conventional approach and can be easily implemented in Matlab/Octave.  The paper is presented in a tutorial way and therefore it include some scripts in Matlab/Octave which demonstrate how easy is the implementation of a Wirtinger's calculus-based formulation of the load flow. The rest of the paper is organized as follows:  Section II presents a brief introduction to complex analysis and Wirtinger's calculus. The proposed load-flow is presented in Section III followed by numerical simulation in Section IV.  Finally, some conclusions (Section V) and relevant references.

\section{Mathematical preliminaries}
\subsection{brief introduction to complex analysis}
Among this paper we represent the complex imaginary unit as $j=\sqrt{-1}$.  Therefore, a complex number can be represented as $z=x+yj$ where $x$ is the real part and $y$ the imaginary part.  A function $f:\mathbb{C}\rightarrow\mathbb{C}$ is a complex-analytical function, also known as  complex differentiable or holomorphic function, if $f(z)$ has a derivative at each point of $\mathcal{B}\subset z$ and if $f(z)$ is single valued.  

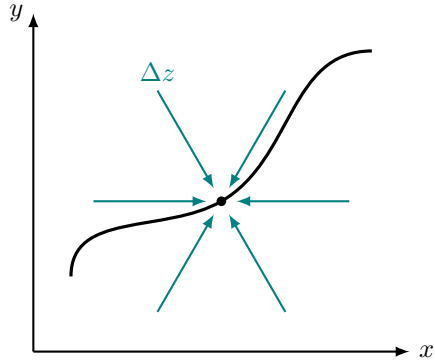
\begin{figure}[htb]
\centering
\begin{tikzpicture}[x=1mm, y=1mm]
\draw[-latex, thick] (-25,-20) -- +(50,0) node[right]{$x$};
\draw[-latex, thick] (-25,-20) -- +(0,45) node[left]{$y$};
\draw[very thick] (0,0) to[out=30,in=180] (20,20);
\draw[very thick] (0,0) to[out=210,in=90] (-20,-10);
\draw[thick,fill] (0,0) circle (0.5);
\draw[thick, blue!50!green,latex-]  (0:2) -- +(0:15);
\draw[thick, blue!50!green,latex-]  (60:2)  -- +(60:15);
\draw[thick, blue!50!green,latex-]  (120:2) -- +(120:15) node[above] {$\Delta z$};
\draw[thick, blue!50!green,latex-]  (180:2) -- +(180:15);
\draw[thick, blue!50!green,latex-]  (-60:2) -- +(-60:15);
\draw[thick, blue!50!green,latex-]  (-120:2) -- +(-120:15);
\end{tikzpicture}
\caption{Different directions to obtain the limit when $\Delta z\rightarrow 0$ in the complex plane}
\label{fig:derivadacompleja}
\end{figure}

The complex derivative can be intuitively defined as 
\begin{equation}
	f'(z) = \lim_{\Delta z \rightarrow 0} \frac{f(z+\Delta z)-f(z)}{\Delta z}
\end{equation}
However, an intermediate problem appears when the function is defined on the complex numbers since there are infinitely many directions in which it is possible to achieve the limit as shown in Fig \ref{fig:derivadacompleja}.   This is similar to the concept of  semi-differentiability in real function where we can obtain a right-derivative different from the left-derivative for discontinuous functions. Therefore, an analytical complex function has the same derivative regardless the direction in which the limit is obtained.
Let us define $z^*=x-yj$ as the complex conjugate of $z$, then a necessary condition for a complex function $f=u(x,y)+v(x,y)$ to be analytic is that 
\begin{equation}
	f'(z^*) = 0
\end{equation}
In other word, the function must depend only on $z$ and not on its conjugate $z^*$.   This condition can be represented in real and imaginary part as follows:
\begin{align}
	\frac{\partial u}{\partial x} &= \frac{\partial v}{\partial y} \\
	\frac{\partial u}{\partial y} &= -\frac{\partial v}{\partial x} 
\end{align}
These are known as the Cauchy-Riemann equations, a set of strong conditions that generates very interesting results on the complex analysis with applications in all sciences including power system engineering \cite{complexbook}. Consequently, the study of complex-analytical function is a mayor area of interest in applied mathematics. However, many practical problems are non-analytical.  In those cases, the conventional approach is to split the problem into real and imaginary part, and use real analysis.

\subsection{Wirtinger's calculus}

\begin{definition}[wirtinger's derivatives]
Given a complex function $f:\mathbb{C}\rightarrow\mathbb{C}$ defined in real and imaginary parts as $f=u+jv$ with $u=u(x,y),v=v(x,y)$, we define the Wirtinger’s derivative, the conjugate Wirtinger’s derivative and the total Wirtinger differential as follows
\begin{align*}
	\frac{\partial f}{\partial z} &= \frac{1}{2}\left( \frac{\partial u}{\partial x}+\frac{\partial v}{\partial y}\right) + \frac{j}{2}\left( \frac{\partial v}{\partial x}-\frac{\partial u}{\partial y}\right) \\
	\frac{\partial f}{\partial z^*} &= \frac{1}{2}\left( \frac{\partial u}{\partial x}-\frac{\partial v}{\partial y}\right) + \frac{j}{2}\left( \frac{\partial v}{\partial x}+\frac{\partial u}{\partial y}\right) \\
	df &= \frac{\partial f}{\partial z}dz + \frac{\partial f}{\partial z^*}dz^*
\end{align*}
\end{definition}

\begin{remark}
When $f$ is complex differentiable then its Wirtinger’s derivative degenerates to the standard complex derivative, while its conjugate Wirtinger’s derivative derivative vanishes\cite{signal}.
\end{remark}

Given a complex function $f=u+jv$ with $u=u(x,y),v=v(x,y)$, we define the Wirtinger’s derivative and the conjugate Wirtinger’s derivative as follows
\begin{align}
	\frac{\partial f}{\partial z} &= \frac{1}{2}\left( \frac{\partial u}{\partial x}+\frac{\partial v}{\partial y}\right) + \frac{j}{2}\left( \frac{\partial v}{\partial x}-\frac{\partial u}{\partial y}\right) \\
	\frac{\partial f}{\partial z^*} &= \frac{1}{2}\left( \frac{\partial u}{\partial x}-\frac{\partial v}{\partial y}\right) + \frac{j}{2}\left( \frac{\partial v}{\partial x}+\frac{\partial u}{\partial y}\right)
\end{align}
When $f$ is complex differentiable then its Wirtinger’s derivative degenerates to the standard complex derivative, while its conjugate Wirtinger’s derivative derivative vanishes\cite{signal}, that is:
\begin{align}
	\frac{\partial f}{\partial z} &= f'(z) \\
		\frac{\partial f}{\partial z^*} &= 0
\end{align}
In general, we are interested in non-analytic functions and therefore $f$ does not satisfied the Cauchy-Riemann equations.  However, we can define easily the following Wirtinger differential
\begin{equation}
	df = \frac{\partial f}{\partial z} + \frac{\partial f}{\partial z^*}
\end{equation}
Each operator behaves similarly as a partial derivative taking into consideration that $f$ depends on both $z$ and its conjugate $z^*$.  Therefore, we can apply common rules for differentiation concerning the sum and the product of functions as follows
\begin{align}
  \frac{\partial(f + g)}{\partial z} &= \frac{\partial f}{\partial z}+\frac{\partial g}{\partial z} \\
	\frac{\partial(f + g)}{\partial z^*} &= \frac{\partial f}{\partial z^*}+\frac{\partial g}{\partial z^*} \\
	\frac{\partial(f\cdot g)}{\partial z} &= f\frac{\partial g}{\partial z}+g\frac{\partial f}{\partial z} \\
	\frac{\partial(f\cdot g)}{\partial z^*} &= f\frac{\partial g}{\partial z}+g\frac{\partial f}{\partial z^*} 
\end{align}
Wirtinger differential also allows to use the chain rule. In addition,  it holds that
\begin{align}
	\frac{\partial z}{\partial z^*} &= 0 \\
	\frac{\partial z^*}{\partial z} &= 0 
\end{align}
for the sake of clarity, let us see a simple example
\begin{align}
	d(|z|^2) &= d(z z^*) \\
			     &= z d(z^*) + z^* d(z) \\
					 &= z\left( \frac{\partial z^*}{\partial z} + \frac{\partial z^*}{\partial z^*} \right) + z^*\left( \frac{\partial z}{\partial z} + \frac{\partial z}{\partial z^*} \right) \\
					&= z (0+1)+z^*(1+0) \\
					&= z+z^*
\end{align}
It is important to remember that $df$ is a Wirtinger differential and not the conventional complex derivative.  In fact, $f=zz^*$ is non-analytical on the complex domain.  The same function can be represented as a real function $f:\mathbb{R}^2\rightarrow \mathbb{R}$ given by $f(x,y) = x^2+y^2$ which is analytical in the real domain.  This is the trick usually applied in load-flow calculations.  The real diferential of $f$ is the gradient given by $\nabla f$, that is
\begin{align}
\nabla f &= \left(\frac{\partial f}{\partial x},\;\frac{\partial f}{\partial y}\right)^T \\
         &= \left(2x,\;2y\right)^T
\end{align}
Notice that, in this case, the Wirtinger differential is exactly equal to the gradient of $f(x,y)$ when the former is spitted in real and imaginary part.

\section{Load-flow in power distribution grids}

Let us consider a power distribution grid represented by its admitance matriz $Y_{BUS}=(y_{km})$.  The slack node is marked with a subscript $0$ and the remaining nodes $N=\left\{1,2,\dots,n\right\}$. All vectors and matrices are represented in uppercase letters whereas their components are represented by lowercase letter.  Voltage in the slack node $v_0$ is known, hence our objective is to obtain the nodal voltages $v_k$ with $k\in N$.  

In the following, all variables and equations are represented in the complex domain and all derivatives are  Wirtinger’s derivatives. For the sake of simplicity we assume initially that all loads admit a constant power model, after that we generalize for the case of ZIP loads.

\subsection{Load-flow with constant power models}

Let us define $v_k$ and $i_k$ as the nodal power in each node $k\in N$, hence, the nodal power is given by
\begin{equation}
	s_k^{*} = v_k^* i_k
\end{equation}
where
\begin{equation}
	i_k = y_{k0}v_0 + y_{kk}v_k + \sum_{m\neq k}^{n} y_{km}v_m
\end{equation}
we use $s_k^*$ instead of $s_k$ for the sake of simplicity. Notice that the function is non-analytic since it depends on both, $v$ and $v^*$. In order to obtain a Newton's step we require to calculate a the following derivatives and conjugate derivatives:
\begin{align}
\frac{\partial s_k^*}{\partial v_k} &= v_k^*\frac{\partial i_k}{\partial v_k} = y_{kk}  v_k^* \\
\frac{\partial s_k^*}{\partial v_k^*} &= i_k \\
\frac{\partial s_k^*}{\partial v_m} &= v_k^*\frac{\partial i_k}{\partial v_m}  = y_{km}  v_k^* \\
\frac{\partial s_k^*}{\partial v_m^*} &= 0
\end{align}
Therefore we can obtain the differential $ds_k^*\approx \Delta s_k$ as
\begin{equation}
	\Delta s_k^* = \frac{\partial s_k^*}{\partial v_k}\Delta v_k + \frac{\partial s_k^*}{\partial v_k^*}\Delta v_k^* + \sum_{k\neq m}^{n}\frac{\partial s_k^*}{\partial v_m}\Delta v_m + \frac{\partial s_k^*}{\partial v_m^*}\Delta v_m^*
\end{equation}
After simple algebraic manipulations we obtain the following matrix representation
\begin{equation}
	\Delta S_N^* = (diag(V_N^*)Y_{NN})\Delta V_N + diag(I_N)\Delta V_N^* \label{eq:Jac1}
\end{equation}
This is clearly a linear function that can be easily solved generating the Newton's iteration described in Algorithm 1.
\begin{algorithm}
\caption{Load-flow based on Wirtinger's calculus}\label{Alg:conpotencia}
\begin{algorithmic}[1]
\Require $YBUS,S_N$ 
\State $V_N \gets 1pu$
\State $I_N \gets Y_{BUS}V_N$
\State $\Delta S_N \gets S_N-V_N\circ I_N^* $ 
\State $\epsilon \gets \left\|\Delta S_N \right\| $
\While{$\epsilon > \textit{tolerance}$}
  \State Solve Eq (\ref{eq:Jac1}) 
	\State $V_N \gets V_N + \Delta V_N$
	\State $\Delta S_N \gets S_N-V_N\circ I_N^* $ 
\EndWhile
\State Calculate $S_{loss}$
\end{algorithmic}
\end{algorithm}
The main advantage of Wirtinger's calculus is that load flow equations can be obtained in a straightforward manner without separate in real and imaginary parts, resulting in a compact representation of the Jacobian given by (\ref{eq:Jac1}) which remains in the complex domain.   Distributed generation can be included in the model by changing the sign of $s_k$. 

\subsection{Solution of the algebraic system}
The most time consuming step in Algorithm 1 is the solution of the linear-complex algebraic system (\ref{eq:Jac1}).  In this aspect, it is possible to apply different techniques in order to obtain a fast algorithm.  First, notice that $|v_k|\neq 0$ (otherwise the grid is in short-circuit), therefore, we can pre-multiply (\ref{eq:Jac1}) by $diag(V_N^*)^{-1}$ obtaining the following equation
\begin{equation}
	\Delta S_N^*/V_N^* = Y_{NN}\Delta V_N + diag(I_N/V_N^*)\Delta V_N^* 
\end{equation}
where $\Delta S_N^*/V_N^*$ and $I_N/V_N^*$ represent element-wise divisions.  Let us define $J_N=\Delta S_N/V_N$ and $K_N=diag(I_N/V_N^*)$ then we can further simplify Equation (\ref{eq:Jac1})
\begin{equation}
	J_N^* = Y_{NN}\Delta V_N + K_N\Delta V_N^* 
\end{equation}
Notice that $J_N$ and $K_N$ is a vector and a diagonal matrix respectively.  Thus, they can be obtained at low computational cost.  By a simple conjugation we can obtain a double-size linear system as follows
\begin{equation}
	\left( \begin{array}{c}J_N^* \\ J_N\end{array}\right) =
	\left( \begin{array}{cc} Y_{NN} & K_N \\ K_N^* & Y_{NN}^*\end{array}\right)
	\left( \begin{array}{c}\Delta V_N \\ \Delta V_N^*\end{array}\right)
	\label{eq:lineal}
\end{equation}
There are many well known techniques for solving systems of this type, specially taking into account that $Y_{NN}$ is highly disperse and can be easily factorized using an LDU technique.

\subsection{Load-flow with ZIP models}

Let us consider now the case of ZIP models.  In this case, the nodal power is given by
\begin{equation}
s_k^* |v_k|^\alpha = v_k^* i_k
\end{equation}
where the real value $\alpha$ indicates the type of model: $0$ for constant power, $1$ for constant current and $2$ for constant impedance. Define now a function $f_k:\mathbb{C}^n \rightarrow \mathbb{C}$ as
\begin{equation}
	f_k = v_k^* i_k - s_k^* (v_kv_k^*)^{\alpha/2} \label{eq:zipmodel}
\end{equation}
with its corresponding Wirtinger derivatives
\begin{align}
\frac{\partial f_k}{\partial v_k} &= y_{kk}  v_k^* - \frac{\alpha}{2}s_k^* |v_k|^{\alpha-2}v_k^* \\
\frac{\partial f_k}{\partial v_k^*} &= i_k -\frac{\alpha}{2}s_k^* |v_k|^{\alpha-2}v_k\\
\frac{\partial f_k}{\partial v_m} &= y_{km}  v_k^* \\
\frac{\partial f_k}{\partial v_m^*} &= 0
\end{align}
hence, we obtain the following differential
\begin{equation}
\begin{split}
	\Delta F_N = diag\left(I_N-\frac{\alpha}{2}\circ S_N^* \circ V_N^{\alpha-2} \circ V_N^*\right)\Delta V_N^* \\
	+\left(diag(V_N^*)Y_{NN}-diag\left(\frac{\alpha}{2}\circ S_N ^*\circ V_N^{\alpha-2} \circ V_N^*\right)\right)\Delta V_N \label{eq:Jac2}
\end{split}
\end{equation}
where $\circ$ is the Hadamard product (i.e the element-wise multiplication in Matlab).  The Newton's iteration can be represented by Algorithm 2.
\begin{algorithm}
\caption{Load-flow based on Wirtinger's calculus with ZIP models}\label{Alg:conZIP}
\begin{algorithmic}[1]
\Require $YBUS,S_N,\alpha$ 
\State $V_N \gets 1pu$
\State $I_N \gets Y_{BUS}V_N$
\State $\Delta F_N \gets $ Eq(\ref{eq:zipmodel})
\State $\epsilon \gets \left\|\Delta F_N \right\| $
\While{$\epsilon > \textit{tolerance}$}
  \State Solve Eq (\ref{eq:Jac2}) 
	\State $V_N \gets V_N + \Delta V_N$
	\State $\Delta S_N \gets S_N-V_N\circ I_N^* $ 
\EndWhile
\State Calculate $S_{loss}$
\end{algorithmic}
\end{algorithm}
Notice that Eqs (\ref{eq:Jac1}) and (\ref{eq:Jac2}) are compact and elegant representations of the power flow which simplify analysis and calculations. The advantages in term of simplicity are evident, however the numerical performance of the algorithm depends highly on the solution of the complex conjugate linear system.  As in the previous case, we can pre-multiply by $diag(V_N^*)^{-1}$ obtaining the following expresion
\begin{equation}
	L_N = (Y_{NN}-H_N)\Delta V_N + (K_N-H_N)\Delta V_N^*
\end{equation}
with
\begin{align}
	L_N &= \Delta F_N/V_N^* \\
	K_N &= diag(I_N/V_N^*) \\
	H_N &= diag(\frac{\alpha}{2} S_N^* V_N^{\alpha-2})
\end{align}
Again, each of these matrices can be calculated at low computational effort with the following linear representation
\begin{equation}
\left(\begin{array}{c}L_N\\L_N^*\end{array} \right)
\left(\begin{array}{cc}Y_{NN}-H_N & K_N-H_N \\K_N^*-H_N^* & Y_{NN}^*-H_N^*\end{array} \right)
\left(\begin{array}{c}\Delta V_N\\\Delta V_N^*\end{array} \right)
\end{equation}
Other representations of this linear system are possible.  More research is required in this area.

\section{Numerical results}

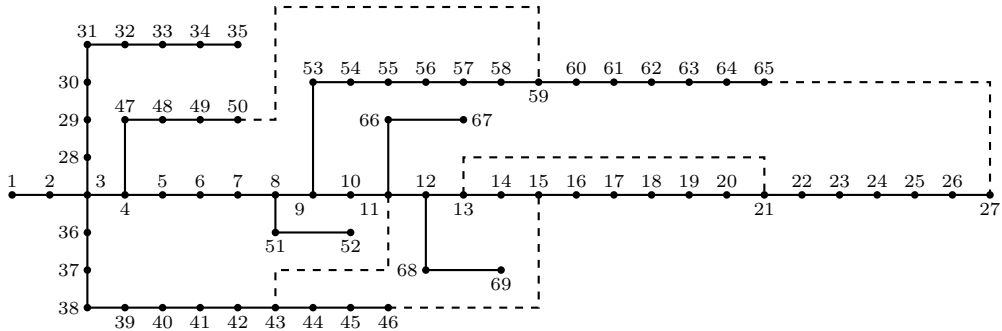
\begin{figure}
\centering
\scriptsize
\begin{tikzpicture}[x=1.0mm,y=1.0mm, thick,rotate=90]
\fill (0,0) circle (0.5) node[above]  {1};
\fill (0,-5) circle (0.5) node[above]  {2};
\fill (0,-10) circle (0.5) node[above right] {3};
 \fill (0,-15) circle (0.5) node[below] {4};
 \fill (0,-20) circle (0.5) node[above] {5};
 \fill (0,-25) circle (0.5) node[above] {6};
 \fill (0,-30) circle (0.5) node[above] {7};
 \fill (0,-35) circle (0.5) node[above] {8};
 \fill (0,-40) circle (0.5) node[below left] {9};
 \fill (0,-45) circle (0.5) node[above] {10};
 \fill (0,-50) circle (0.5) node[below left] {11};
 \fill (0,-55) circle (0.5) node[above] {12};
 \fill (0,-60) circle (0.5) node[below] {13};
 \fill (0,-65) circle (0.5) node[above] {14};
 \fill (0,-70) circle (0.5) node[above] {15};
 \fill (0,-75) circle (0.5) node[above] {16};
 \fill (0,-80) circle (0.5) node[above] {17};
 \fill (0,-85) circle (0.5) node[above] {18};
 \fill (0,-90) circle (0.5) node[above] {19};
 \fill (0,-95) circle (0.5) node[above] {20};
 \fill (0,-100) circle (0.5) node[below] {21};
 \fill (0,-105) circle (0.5) node[above] {22};
 \fill (0,-110) circle (0.5) node[above] {23};
 \fill (0,-115) circle (0.5) node[above] {24};
 \fill (0,-120) circle (0.5) node[above] {25};
 \fill (0,-125) circle (0.5) node[above] {26};
 \fill (0,-130) circle (0.5) node[below] {27};
 \fill (5,-10) circle (0.5) node[left] {28};
 \fill (10,-10) circle (0.5) node[left] {29};
 \fill (15,-10) circle (0.5) node[left] {30};
 \fill (20,-10) circle (0.5) node[above] {31};
 \fill (20,-15) circle (0.5) node[above] {32};
 \fill (20,-20) circle (0.5) node[above] {33};
 \fill (20,-25) circle (0.5) node[above] {34};
 \fill (20,-30) circle (0.5) node[above] {35};
 \fill (-5,-10) circle (0.5) node[left] {36};
 \fill (-10,-10) circle (0.5) node[left] {37};
 \fill (-15,-10) circle (0.5) node[left] {38};
 \fill (-15,-15) circle (0.5) node[below] {39};
 \fill (-15,-20) circle (0.5) node[below] {40};
 \fill (-15,-25) circle (0.5) node[below] {41};
 \fill (-15,-30) circle (0.5) node[below] {42};
 \fill (-15,-35) circle (0.5) node[below] {43};
 \fill (-15,-40) circle (0.5) node[below] {44};
 \fill (-15,-45) circle (0.5) node[below] {45};
 \fill (-15,-50) circle (0.5) node[below] {46};
 \fill (10,-15) circle (0.5) node[above] {47};
\fill (10,-20) circle (0.5) node[above] {48};
\fill (10,-25) circle (0.5) node[above] {49};
\fill (10,-30) circle (0.5) node[above] {50};
\fill (-5,-35) circle (0.5) node[below] {51};
\fill (-5,-45) circle (0.5) node[below] {52};
\fill (15,-40) circle (0.5) node[above] {53};
\fill (15,-45) circle (0.5) node[above] {54};
\fill (15,-50) circle (0.5) node[above] {55};
\fill (15,-55) circle (0.5) node[above] {56};
\fill (15,-60) circle (0.5) node[above] {57};
\fill (15,-65) circle (0.5) node[above] {58};
\fill (15,-70) circle (0.5) node[below] {59};
\fill (15,-75) circle (0.5) node[above] {60};
\fill (15,-80) circle (0.5) node[above] {61};
\fill (15,-85) circle (0.5) node[above] {62};
\fill (15,-90) circle (0.5) node[above] {63};
\fill (15,-95) circle (0.5) node[above] {64};
\fill (15,-100) circle (0.5) node[above] {65};
\fill (10,-50) circle (0.5) node[left] {66};
\fill (10,-60) circle (0.5) node[right] {67};
\fill (-10,-55) circle (0.5) node[left] {68};
\fill (-10,-65) circle (0.5) node[below] {69};

\draw (0,0) -- (0,-130);
\draw (0,-10) -| (20,-30);
\draw (0,-10) -| (-15,-50);
\draw (0,-15) -| (10,-30);
\draw (0,-35) -| (-5,-45);
\draw (0,-40) -| (15,-100);
\draw (0,-50) -| (10,-60);
\draw (0,-55) -| (-10,-65);
\draw[dashed] (10,-30) |- (25,-35) |- (15,-70);
\draw[dashed] (0,-50) -| (-10,-35) |- (-15,-35);
\draw[dashed] (-15,-50) |- (0,-70);
\draw[dashed] (0,-60) -| (5,-100) -- (0,-100);
\draw[dashed] (15,-100) |- (0,-130);
\end{tikzpicture}
\caption{IEEE 69-node test system}
\label{fig:testsystem}
\end{figure}

The proposed algorithms were evaluated on the IEEE 69-nodes test system \cite{testsystem69} depicted in Fig \ref{fig:testsystem}.  This test system was designed for network reconfiguration and therefore it has some tie lines. Initially, the load-flow was calculated considering all nodes as constant power loads and tie lines open. Algorithm 1 converged in only  3 iterations with a tolerance of $1\times 10^{-4}$. The script was implemented in Matlab2016 and computation time was $0.012059$ seconds in a computer with processor Intel(R) Core(TM) i7-6700CPU@3.40GHz with RAM8.00Gb and windows 64bits operative system. 

A second numerical simulation was performed but this time the tie lines were connected resulting a meshed power distribution grid.  Algorithm 1converged in 3 iterations with the same tolerance.  Computation time was $0.015036$ seconds.  It is important to notice the proposed algorithm does not require any modification for meshed grids.

In order to maintain the tutorial presentation of this paper, we include the code for Algorithm 1 in Matlab (see appendix A). Notice the simple implementation of the algorithm thanks to use of the capability of matrix calculations in the complex domain.

A third set of numerical calculations was performed but this time, the grid was considered with ZIP loads and renewable energy generation.   Convergence was achieved in 3 iterations and computation time was 0.007484 seconds. Matlab code is available in Appendix B.
\section{Conclusions}

A general methodology for the load-flow calculations in power distribution grids was presented.  This methodology uses Wirtinger's calculus instead of conventional complex analysis allowing a compact representation of the algebraic equations and a simple implementation in scripting based languages such as Matlab/Octave. 

The proposed method allows to include constant power loads, ZIP models and distributed generation. In addition, meshed distribution grids can be directly implemented without any change in the code or the convergence performance.  Numerical results in Matlab demonstrated the methodology is easily implementable and allows fast solutions.

It is the author belief, this approach is very useful from the pedagogical point of view, since the equations can be effortless obtained and interpreted.  Moreover, the implementation of the method is intuitive and this simplify debugging process, especially in complex algorithms were the load flow is just one sub-routine.   The use of Wirtinger's calculus could be extended to other areas of power system analysis, optimization and control.  Possible applications include economical dispatch, optimal power flow, stability and short-circuit analysis.

\bibliographystyle{elsarticle-num-names}
\bibliography{bibliografia}

\appendix
\subsection{Script in Matlab for Algorithm 1}

\begin{lstlisting}[frame=trBL]
tic
Vn = ones(n-1,1)*Vo;
In = Yno*Vo + Ynn*Vn;
dS = Sn - Vn.*conj(In); 
epsilon = norm(dS);
while epsilon>1E-4    
   Jn = dS./Vn;      
   Kn = diag(In./conj(Vn));
   DV = linsolve([Ynn,Kn;Kn',Ynn'],[conj(Jn);Jn]);
   Vn = Vn + DV(1:n-1);
   In = Yno*Vo + Ynn*Vn;
   dS = Sn - Vn.*conj(In); 
   epsilon = norm(dS);
end
toc
\end{lstlisting}

\subsection{Script in Matlab for Algorithm 2}

\begin{lstlisting}[frame=trBL]
tic
Vn = ones(n-1,1)*Vo;
In = Yno*Vo + Ynn*Vn;
dFn = conj(Vn).*In-conj(Sn).*abs(Vn).^alpha;
epsilon = norm(dFn);
while epsilon>1E-4
   Ln = dFn./conj(Vn);
   Kn = diag(In./conj(Vn));
   Hn = diag(alfa/2.*conj(Sn).*Vn.^(alpha-2));       
   DV = linsolve([Ynn-Hn,Kn-Hn;Kn'-Hn',Ynn'-Hn'],...
		     [Ln;conj(Ln)]);
   Vn = Vn - DV(1:NumN-1);
   In = Yno*Vo + Ynn*Vn;
   dFn = conj(Vn).*In-conj(Sn).*abs(Vn).^alpha;    
   epsilon = norm(dFn);
end
toc
\end{lstlisting}
\end{document}